\documentclass[11pt]{amsart}

\usepackage{amssymb}
\usepackage{graphicx}
\usepackage{caption} 
\usepackage{subcaption}
\usepackage{amsmath}
\usepackage{amsfonts}
\usepackage{amscd}
\usepackage{color}
\usepackage{verbatim}
\usepackage{amsthm}
\usepackage{marginnote}
\usepackage{enumerate}


\theoremstyle{definition}

\numberwithin{equation}{section}

\begin{document}
\date{\today}
\title[Dynamics on immediate basins for parabolic Newton maps]{ Dynamics on immediate basins for parabolic Newton maps}

\author[K. Mamayusupov]{Khudoyor Mamayusupov}
\address{National Research University Higher School of Economics, Russian Federation \\
Faculty of Mathematics \\
 Ul. Usacheva 6, Moscow, Russian Federation }
\email{kmamayusupov@hse.ru}
\thanks{Research was partially supported by the ERC advanced grant ``HOLOGRAM'' and the Deutsche Forschungsgemeinschaft SCHL 670/4}
\begin{abstract} Dynamics on parabolic immediate basins for rational Newton maps of entire functions have been studied. It is proved that every parabolic immediate basin contains invariant accesses to the parabolic fixed point at infinity. Moreover, among these accesses there exists a unique dynamically defined access where dynamics are attracted towards the parabolic fixed point, whereas for other accesses, if there is any, the dynamics are repelled.  
\end{abstract}
\maketitle
\section{Introduction}

Let $f$ be a complex polynomial or a transcendental entire function defined on the complex plane $\mathbb C$. A meromorphic function $N_f(z):=z-f(z)/f'(z)$ is called the Newton map of $f(z)$. In the literature, the Newton map of a polynomial $p(z)$ is also called the Newton method applied to $p(z)$. For a polynomial equation $p(z)=0$, if one considers the corresponding Newton map $N_p(z)$ then the roots of $p(z)$ in $\mathbb C$ become the attracting fixed points of its Newton map $N_p(z)$. Hence searching for the roots of a complex polynomial is equivalent to searching for the attracting fixed points of the corresponding Newton map. At a starting point $z_0$ in the complex plane $\mathbb{C}$, we iteratively apply the Newton map $N_f(z)$ to produce the sequence \[z_0, N_f(z_0), N_f^{\circ 2}(z_0), \dots, N_f^{\circ n}(z_0), \dots
  \]
Here we write $F^{\circ n}$ for the $n$th iterate of a complex map $F$, for instance, $F^{\circ 2}(z):=F(F(z))$. 
If this sequence converges to a complex number $\xi$ then $\xi$ is a fixed point of $N_f$ and moreover the same sequence converges to a root of $f(z)=0$ and approximates the simple roots in quadratic order, thus the convergence is very rapid. In this paper, 
we are interested in the case when the Newton map is rational. We do not study here the numerical aspects of Newton's method (refer to the works [9, 10]) but rather consider the Newton maps as dynamical systems on the Riemann sphere $\hat{\mathbb C}$. More precisely, we study dynamical properties of rational Newton maps $N_f$ for $f=p e^q$ with non-constant $q$ on the open set of the complex plane where the above defined sequence converges in $\hat{\mathbb C}$. This set is called the basin of attraction of a fixed point $\xi$ of $N_f$. It is an open set and belongs to the stable set called the Fatou set of $N_f$, the complement of which is called the Julia set where the dynamics is chaotic (sensitive to initial conditions). The Fatou set of a holomorphic function $F$ is defined as the set of $z\in \hat{\mathbb C}$ such that there exists $U$ an open neighborhood of $z$ where the iterates $\{F^{\circ n}|_U, n\ge 1\}$ form a normal family on $U$. The Fatou set of a rational function $F$ coincides with the set of points $z\in \hat{\mathbb C}$ {\em stable in the sense of Lyapunov} for $F$.  

Let us introduce some basic definitions. For a fixed point $F(z)=z$ the quantity $\lambda=F'(z)$ is a local conformal invariant (under a conformal conjugacy) and is called the {\em multiplier} of $z$. The fixed point $z$ is called {\it attracting} if $|\lambda|<1$, in particular, {\it superattracting} if $\lambda=0$, {\it repelling} if $|\lambda|>1$, {\it indifferent} if $|\lambda|=1$, in this case let $\lambda=e^{2\pi i \theta}$, then {\it rationally indifferent} (also called parabolic) if $\theta\in\mathbb{Q}$, {\it irrationally indifferent} if $\theta\not\in\mathbb{Q}$.
The following is known for fixed points. Attracting fixed points and irrationally indifferent fixed points if the function locally linearisable (Siegel point) belong to the Fatou set. Repelling and parabolic fixed points belong to the Julia set. Irrationally indifferent fixed points if the function is not locally linearisable (Cremer point) belong to the Julia set as well.
The multiplier is also defined for periodic points and the same classification and the similar properties are true for them as well.
The multiplicity of a fixed point $z_0$ is defined as the multiplicity of $z_0$ as a root of the fixed point equation $F(z)=z$. If $z=0$ is a parabolic fixed point with the multiplier $+1$ then $F(z)=z+a z^{m+1}+O(z^{m+2})$, where $a\neq 0$ is a complex number, is a Taylor series expansion of $F$ near the origin and $m+1\ge 2$ is the multiplicity of the parabolic fixed point at the origin. The number $m$ in this case is called the parabolic multiplicity of a parabolic fixed point.  

Every rational function of degree at least $2$ has a fixed point that is either repelling or parabolic with the multiplier $+1$. We call this type of fixed points \emph{weakly repelling}. If such a point is unique then the Julia set is connected by theorem of Shishikura [11]. For Newton maps the point at $\infty$ is the only weekly repelling fixed point, as a corollary we obtain that the Julia set for all rational Newton maps is connected.

Denote $\deg(F,z)$ the local degree of a function $F$ at a point $z$. Denote $C_F =\{z|\deg(F,z)>1 \}$ and $P_F=\overline{\bigcup_{n\ge1} F^{\circ n}(C_F)}$ the set of critical points and the post-critical set of $F$ respectively. For holomorphic functions, the set of critical points is exactly the set where the derivative vanishes. A function $F$ is called \emph{post-critically finite} if $P_F$ is finite. A function $F$ is called \emph{geometrically finite} if the intersection of $P_F$ with the Julia set is a finite set. 

For geometrically finite rational Newton maps the Julia sets are locally connected, it is corollary of a result in [12], where it asserts that for geometrically finite rational functions the Julia set is locally connected if it is connected.

The rational Newton maps can be described in terms of multiplies of fixed points as in the following theorem in [8].

\textbf{Theorem 1.} \emph{Let $N:\mathbb C \to \hat{\mathbb C}$ be a meromorphic function. It is the Newton map of an entire function $f:\mathbb C\to\mathbb C$ if and only if for each fixed point $\xi$, $N(\xi)=\xi$, there is an integer $m=m_{\xi} \ge 1$ such that $N'(\xi)=(m-1)/m$. In this case there exists a constant  $c\in \mathbb C \backslash \{0\}$ such that $f=c \cdot e^{\int\frac{1}{\zeta-N(\zeta)}d \zeta }$. Entire functions $f$ and $g$ have
the same Newton map if and only if $f=c\cdot g$ for some constant $c\in \mathbb C \backslash \{0\}$.
}

The point at $\infty$ is a removable singularity for some Newton maps making them rational functions on $\hat{\mathbb C}$ [8].

\textbf{Theorem 2.} \emph{ Let $f:\mathbb C\to\mathbb C$ be an entire function. Its Newton map $N_f$ is a rational function if and only if there are polynomials $p$ and $q$ such that $f$ has the form $f=p e^q$. 
Let $m$ and $n$ be the degrees of $p$ and $q$, respectively. For $n = 0$ and $m \geq 2$, the point at $\infty$ is repelling with the multiplier $m/(m-1)$. For the pair $n = 0$ and $m = 1$, $N_f$ is constant. For $n > 0$, the point at $\infty$ is parabolic with the multiplier $+1$ and the multiplicity $n+1 \geq 2$.
}

\section{Results}

 Here is a criterion to check whether or not a given rational map is a Newton map. It is based on partial fraction decomposition of rational functions.

\textbf{Theorem 3.} [Description of rational Newton maps]\emph{
Let a rational function $N :  \hat{\mathbb C}\to \hat{\mathbb C}$ of degree $d\geq 2$ be given. Assume $\infty$ is a weakly repelling fixed point of $N$. Let a partial fraction decomposition: $\cfrac{1}{z-N(z)}=\sum_{i=1}^k r_i(\cfrac{1}{z-z_i})+s(z)$ be given, where $r_i=r_i(z)$, for $1\leq i\leq k$, and $s=s(z)$ are polynomials with normalizations $r_i(0)=0$ for $1\leq i\leq k$, where $z_i$ runs over all distinct fixed points of $N$ in $\mathbb{C}$. Then $N$ is a Newton map of an entire function if and only if there exist integers $m_i\ge 1$ such that $r_i(z)\equiv m_i\cdot z$. In this case, let $p = p(z) =\prod_{i=1}^k (z-z_i)^{m_i}$ (if there is no finite fixed point of N then we let $p=1$) and $q=q(z) := \int_0^z s(w)dw$ be polynomials, then $N = N_{p e^q}$ and its degree $d=k+\deg(q)$.
}

\begin{proof}
	Let $N :  \hat{\mathbb C}\to \hat{\mathbb C}$  be a rational function of degree $d\geq 2$, without loss of generality assume that $ \infty$ is a weekly repelling fixed point (we can always conjugate $N$ by a suitable M\"obius map sending its weekly repelling fixed point to $\infty$). Assume the partial fraction decomposition of $\cfrac{1}{z-N(z)}$ is given $$\frac{1}{z-N(z)}=\sum_{i=1}^k r_i(\frac{1}{z-z_i})+s(z)$$ and assume $r_i(z)\equiv m_i\cdot z$ for integers $m_i\ge 1$, (see for details Chapter~1.4 in [1]). We have $$\cfrac{1}{z-N(z)}=\sum_{i=1}^k \frac{m_i}{z-z_i}+s(z).$$ If we denote polynomials $p=p(z):=\prod_{i=1}^k (z-z_i)^{m_i}$ and $q = q(z) := \int_0^z s(w)dw$, the latter reduces to $N(z)=z-p(z)/(p'(z)+p(z)\cdot q'(z))$. It follows that $N$ is the Newton map of the entire function $p e^q$ by the uniqueness of Newton maps (Theorem~1). 
	
	Conversely, let $N$ be the Newton map of an entire function $f=p e^q$ for polynomials $p$ and $q$. Let $p(z):=\prod_{i=1}^k (z-z_i)^{m_i}$, where $z_i$ runs over all distinct roots of $p$, then we obtain $$\frac{1}{z-N_f(z)}=\frac{f'}{f}=\frac{p'e^q+pq'e^q}{pe^q}=\frac{p'+pq'}{p}=\frac{p'}{p}+q'=\sum_{i=1}^k \frac{m_i}{z-z_i}+q'(z).$$
	The result follows by uniqueness of partial fraction decompositions for rational functions.
	
	It remains to show that the degree of the Newton map is equal the number of distinct roots of $p$ plus the degree of $q$.
	 We have $$N_{p e^q}(z)=z-\frac{p(z)}{p'(z)+p(z)\cdot q'(z)}=\frac{z(p'(z)+p(z)\cdot q'(z))-p(z)}{p'(z)+p(z)\cdot q'(z)}.$$ 

If the numerator and denominator of the latter has some cancellation factor, then the polynomial equations $z(p'(z)+p(z) q'(z))-p(z)=0$ and $p'(z)+p(z) q'(z)=0$ have a common solution for some $z=z_0$.
	
Plugging the second into the first, we get $p(z_0)=0$. Combining it with the second, we derive $p'(z_0)=0$. These mean that $z=z_0$ is a multiple root of $p$. Thus, the degree $d$ of the Newton map equals $k+deg(q)$, where $k$ is the number of distinct roots of $p$.
\end{proof}
The following notion is the main object when one studies Newton maps. Since we are dealing with a general family of rational Newton maps, we allow a parabolic fixed point at $\infty$. 

\textbf{Definition 1.} [The Basin of Attraction] Let $\xi$ be an attracting or a parabolic fixed point of $F$. The \emph{basin} $\mathcal{A}(\xi)$ of $\xi$ is defined as an interior of the set of starting points that converge to $\xi$ under the iterates of $F$, it is $\text{int}\{z \in \hat{\mathbb{C}}:\lim_{n\rightarrow \infty}F^{\circ n}(z)=\xi\}$. The \emph{immediate basin}  of $\xi$, denoted $\mathcal{A}^{\circ}(\xi)$, is the forward invariant connected component of the basin. 

\begin{figure}[ht!]
	\centering
	\includegraphics[scale=.3]{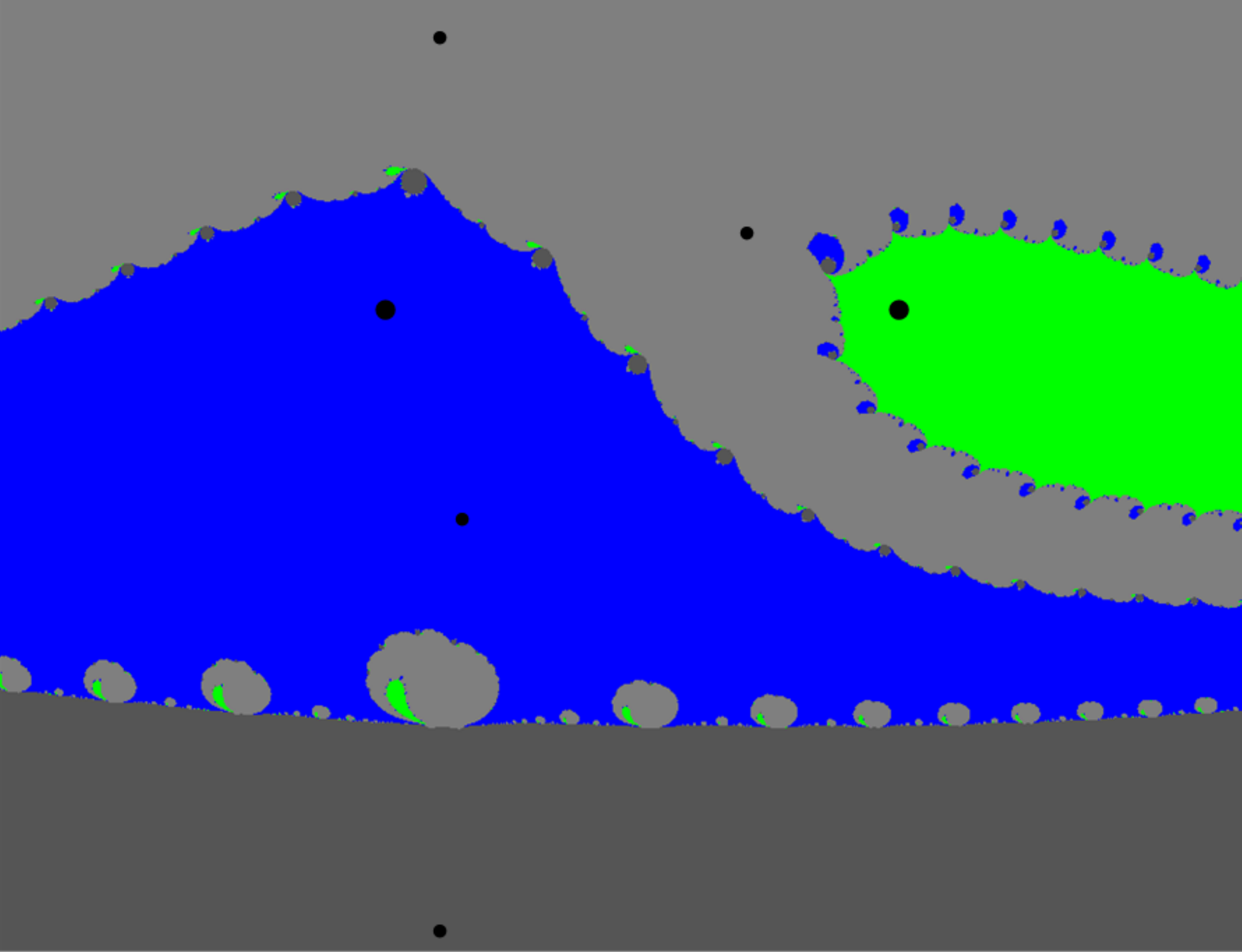}
	\caption{The Fatou components of a quartic Newton map. Blue and green regions are attracting basins, gray regions are the basin of infinity with two distinct immediate basins. Critical points are drawn as black dots with the two heavy dots representing the two superattracting fixed points. There are $2$ critical points in the blue basin and in the gray basin on the top with corresponding $2$ accesses to $\infty$ in each.}
\end{figure}

\textbf{Definition 2.} [Invariant access to $\infty$] Let $\mathcal{A}^{\circ}$ be the immediate basin of a fixed point $\xi \in \mathbb{C}$ or the parabolic fixed point at $\infty$ of a Newton map $N_{p e^q}$. Fix a point $z_0\in \mathcal{A}^{\circ}$, and consider a curve $\Gamma:[0,\infty)\to \mathcal{A}^{\circ}$ with $\Gamma(0) =z_0$ and $\lim_{t\to\infty}\Gamma(t)=\infty$. Its homotopy class (with endpoints fixed) within
$\mathcal{A}^{\circ}$ defines an {\em access to $\infty$} in $\mathcal{A}^{\circ}$. An access is called invariant if together with every representative $\Gamma$ of the access the image $N_{p e^q}(\Gamma)$ also belongs to the same access. 

In the case of an attracting immediate basin $\mathcal{A}^{\circ}(\xi)$, the end point $z_0$ of the homotopy is at the attracting fixed point $\xi$. 

\textbf{Remark 1.} 
	For parabolic immediate basins the invariant accesses are well defined. Indeed, if we choose a point $z_0$ as one of the end points of the homotopy (the other end is at $\infty$), then within the immediate basin by considering a composition of some curve joining $z_0$ and $N_{p e^q}(z_0)$ and the image curve $N_{p e^q}(\Gamma)$ we are always in the same homotopy class of the curve $\Gamma$ thanks to the simply connectivity of the immediate basins.

The invariant accesses are defined for immediate basins of attracting fixed points for Newton maps of polynomials in [3]. Our definition of an invariant access for these domains coincides with theirs.

The following main theorem gives the structure of immediate basins of rational Newton maps. This result is proved in full generality for meromorphic Newton maps in [2] with different methods than ours. Here, we show that all invariant access come in one-to-one correspondence for each fixed point of a corresponding proper map of the unit disk. Moreover, for parabolic immediate basins there always exists a distinguished access called the {\em dynamical invariant access to $\infty$}. To obtain this access in the parabolic immediate basin, we take any curve $\eta$, starting at some $z_0$ and ending at $N_{p e^q}(z_0)$, and considering the homotopy class of the curve $\Gamma:=\bigcup_{n\ge 0}N_{p e^q}^{\circ n}(\eta)$. This curve lands at $\infty$ and it is forward invariant under the Newton map $N_{p e^q}$ by its definition. This notion of dynamical invariant access is used as a main ingredient to construct a natural bijection between distinct spaces of postcritically finite and postcritically minimal Newton maps using parabolic surgery in [4,5,6]. Both types of Newton maps are structurally important considered in their corresponding parameter planes. 

\textbf{Theorem 4.} [Invariant accesses to $\infty$]  \emph{ Let $N_{p e^q}$ be a rational Newton map of degree $d\geq 2$ and $\mathcal{A}^{\circ}$ be an immediate basin of a fixed point $\xi$ of $N_{p e^q}$ (attracting or parabolic). If $\mathcal{A}^{\circ}$ contains $k$ critical points of $N_{p e^q}$ counting with multiplicities, then
$N_{p e^q}|_{\mathcal{A}^{\circ}}$ is a branched covering map of degree $k+1$, and $\mathcal{A}^{\circ}$ has exactly $k$ invariant accesses to $\infty$. If $\mathcal{A}^{\circ}$ is a parabolic immediate basin then among its invariant accesses to $\infty$ there always exists a unique dynamical access.
}

\begin{proof} Let $N_{p e^q}$ a rational Newton map and $\mathcal{A}^{\circ}$ its immediate basin be given. There are two cases; case of basins for attracting fixed points or case of basins for the parabolic fixed point at $\infty$. The first case for Newton maps of polynomials was treated in Proposition~6 in [3]. Since arguments in the proof use only local dynamics of the function within the basin their result is also true for attracting immediate basins of rational Newton maps $N_{p e^q}$ with non-constant $q$. 
	
It remains to prove the theorem for parabolic immediate basins of rational Newton maps $N_{p e^q}$. It is the case if $\deg(q)>0$. Following [3], denote $\mathbb{D}=\{z,|z|<1\}$ the unit disk and its boundary $\mathbb{S}^1=\{z,|z|=1 \}$ the unit circle, and let $\mathcal{A}^{\circ}$ be one of the immediate basins of $\infty$. Consider the Riemann map $\psi: \mathcal{A}^{\circ}\to \mathbb{D}$ uniquely determined with $\psi(c)=0$ and $\psi'(c)>0$, where $c$ is any point in $\mathcal{A}^{\circ}$. Then the composition $f=\psi\circ N_{p e^q} \circ \psi^{-1}$ is a proper map of the unit disc $\mathbb{D}$ with the degree, which is equal to the degree of the restriction $N_{p e^q}|_{\mathcal{A}^{\circ}}$. The critical points of $N_{p e^q}$ in $\mathcal{A}^{\circ}$ are mapped to the critical points of $f$ in $\mathbb{D}$ preserving multiplicities. Assume that $N_{p e^q}$ has $k\ge 1$ critical points in $\mathcal{A}^{\circ}$ counting with multiplicities, there is at least one critical point in every immediate basin. The connectivity of the Julia set implies that $\mathcal{A}^{\circ}$ is simply connected, then by the Riemann-Hurwitz formula the degree of the restriction $N_{p e^q}|_{\mathcal{A}^{\circ}}$  is $k+1$. 
	
The map $f$, as a proper self map of $\mathbb{D}$, is a Blaschke product (the product of a finite number of conformal automorphisms of $\mathbb{D}$, see [7]), hence it has an extension to $\hat{\mathbb{C}}$, denote the extension again by $f$. Both $f$ and the restriction $N_{p e^q}|_{\mathcal{A}^{\circ}}$ have the same degree. Then $f$ has $k+2$ fixed points, one of which is a double parabolic (of multiplicity $3$ or of parabolic multiplicity $2$), since we have a parabolic dynamics in $\mathbb{D}$, and the other distinct $k-1$ fixed points are simple and repelling with real multipliers, and all are located on the unit circle. The unit disk $\mathbb{D}$, the unit circle $\mathbb{S}^1$ and $\hat{\mathbb{C}}\setminus \bar{\mathbb{D}}$ are invariant by $f$. Since $f$ can not have critical point on $\mathbb{S}^1$, it is a covering map of $\mathbb{S}^1$ of degree $k+1$, and the orbit for every $z\in \hat{\mathbb{C}}\setminus \mathbb{S}^1$ converges to the unique parabolic fixed point on $\mathbb{S}^1$. Thus the Julia set is the unit circle $\mathbb{S}^1$. Alternatively, observe that the Fatou set of the Blaschke product is invariant under the suitable involution. If the Julia set is not connected then the Fatou set of this Blaschke product has a unique component which includes a part of the unit circle $\mathbb{S}^1$ and the complement of $\mathbb{S}^1$. But then on the Fatou points in $\mathbb{S}^1$ the iterates must eventually converge to the parabolic fixed point on $\mathbb{S}^1$ then necessarily they fall on the repelling petal and gets pushed back from the parabolic fixed point, which leads to a contradiction to the convergence.

As it was mentioned in Remark~1, for accesses it suffices to consider locally at $\infty$ a part of homotopies fixing one end point at the point at $\infty$ in $\mathcal{A}^{\circ}$. The Riemann map $\psi: \mathcal{A}^{\circ}\to \mathbb{D}$ transports homotopies to the unit disk $\mathbb{D}$. The linearizing coordinates of $k-1$ repelling fixed points of $f$ define $k-1$ invariant accesses among $k$ invariant accesses, and the other invariant access comes from a Fatou coordinate of the parabolic fixed point on $\mathbb{S}^1$. The invariant access associated to the parabolic fixed point defines the \emph{dynamical} access. 
	
Assume that the boundary of $\mathcal{A}^{\circ}$ is locally connected, which is true if the Newton map is geometrically finite. Carath\'eodory theorem assures that $\psi^{-1}$ the inverse map to $\psi: \mathcal{A}^{\circ}\to \mathbb{D}$ extends to the closed unit disk $\overline{\mathbb{D}}$ as a continuous map. Denote the extension again by $\psi^{-1}$. By continuity, we obtain $\psi^{-1}\circ f=N_{p e^q}\circ\psi^{-1}$ on $\overline{\mathbb{D}}$. Counted with multiplicities, $k+2$ fixed points of $f$ correspond to $k+2$ fixed points of $N_{p e^q}$ on $\partial\mathcal{A}^{\circ}$. A fixed point of $N_{p e^q}$ that is on the boundary of an immediate basin is the only parabolic fixed point at $\infty$, so the domain $\mathcal{A}^{\circ}$ has invariant accesses to $\infty$ through $k$ distinct directions ($k-1$ directions for the $k-1$ simple and one direction for the triple fixed point of $f$). 

In the case when $\mathcal{A}^{\circ}$ is not locally connected, so that the inverse to the Riemann map does not extend continuously to the closed unit disk, the statement still holds true. Consider a Koenigs coordinate of a repelling fixed point $\xi_j$ that conjugates $f$ locally near the point  $\xi_j$ to the linear map $z\mapsto f'(\xi_j)z$, we take a segment of a straight-line through the origin, which is invariant. We take an invariant curve in the petal associated to the parabolic fixed point of $f$. Let $\gamma$ be the preimage of this curve that lands at $\xi_j$ in the dynamical plane of $f$. Then we have $\gamma\subset f(\gamma)$. Now we pull the curve $\gamma$ by the Riemann map $\psi$ to $\mathcal{A}^{\circ}$. The accumulation set of $\psi^{-1}(\gamma)$ in $\partial\mathcal{A}^{\circ}$ is connected (see Section~$17$ in [7]) and since $\gamma$ is invariant we conclude that the accumulation set is pointwise fixed by $N_{p e^q}$. But $\infty$ is the only fixed point on the Julia set. This gives us $k$ invariant accesses to $\infty$ in $\mathcal{A}^{\circ}$. We need to show that they are all distinct and the only ones.

It is clear that simple curves within $\mathbb{D}$ converging to a given fixed point of $f$ are homotopic so that every fixed point of $f$ defines a unique access in $\mathcal{A}^{\circ}$. Different fixed points of $f$ lead to non-homotopic curves in $\mathcal{A}^{\circ}$ and thus to different accesses. Indeed, let $l_i$,~$l_j\subset \mathbb{D}$ be the radial lines converging to $\xi_i\not =\xi_j$ respectively, parametrized by the radius. Assume by contrary that $\psi^{-1}(l_i)$ and $\psi^{-1}(l_j)$ are homotopic curves in $\mathcal{A}^{\circ}$ by a homotopy fixing end points; $\psi^{-1}(l_i(1))=\psi^{-1}(l_j(1))=\infty$, then one of the components bounded by a simple closed curve $\psi^{-1}(l_i)\cup\psi^{-1}(l_j)$ must be contained in $\mathcal{A}^{\circ}$. Call this component $V$; then $\psi(V)$ must be one of the sectors bounded by $l_i$ and $l_j$; call it $S$. Both $V$ and $S$ are Jordan domains, so $\psi^{-1}$ extends as a homeomorphism from $\bar S$ onto $\bar{V}$ by Carath\'eodory theorem; but then the extension sends the set $\mathbb{S}^1\cap S$ nowhere.

Conversely, we show that every invariant access to $\infty$ in $\mathcal{A}^{\circ}$ comes from a fixed point of the corresponding $f$. Let $\Gamma:[0,1]\to \mathcal{A}^{\circ}\cup\infty$ be a curve representing an access. Then $\psi(\Gamma)$ lands at a point $\upsilon\in \mathbb{S}^1$ by Corollary~$17.10$ in [7], define it as the associated point of $\Gamma$. Then for every $n\ge 1$, $N^{\circ n}_{p e^q}(\Gamma)$ represents an access and thus has its associated point $\upsilon_n\in \mathbb{S}^1$. Since the Newton map $N_{p e^q}$ has a parabolic fixed point at $\infty$, it is locally a homeomorphism there and every fixed point of $f$ gives rise to an access, all $\upsilon_n$ must be contained in the same connected component of $\mathbb{S}^1$ with the fixed points removed; this component is an interval, say $I$, on which $\{\upsilon_n\}$ must be a monotone sequence converging under $f$ to a fixed point $\upsilon$ of $f$ in $\bar{I}$, i.e. to one of the endpoints. If $\upsilon$ is a one of the repelling fixed points of $f$ then it is impossible. For the remaining possibility, assume $\upsilon$ is a parabolic fixed point of $f$ then the sequence $\{\upsilon_n\}\subset \mathbb{S}^1$ converges to the parabolic fixed point, which is also not possible since $\mathbb{S}^1$ is the Julia set of $f$ and every orbit that converges to a parabolic fixed point must follow the attracting direction which is perpendicular to $\mathbb{S}^1$.  

Finally, observe that the unique invariant access corresponding to the parabolic fixed point of $f$ gives rise to the dynamical invariant access and it always exists and unique.
\end{proof}
\subsection*{Acknowledgements}
The author thanks Dierk Schleicher and Vladlen Timorin for their comments that helped to improve the paper. Research was partially supported by the ERC advanced grant ``HOLOGRAM'' and the Deutsche Forschungsgemeinschaft SCHL 670/4.

\end{document}